\documentclass[%
reprint,
groupedaddress,
amsmath,amssymb,
aps,
pre,
floatfix,
longbibliography,
raggedbottom
]{revtex4-1}

\usepackage{graphicx}
\usepackage{dcolumn}
\usepackage{bm}
\usepackage{color}

\usepackage{hyperref}

\begin{document}


\title{Cutting and shuffling with diffusion: Evidence for cut-offs in interval exchange maps}

\author{Mengying Wang}%
\affiliation{School of Mechanical Engineering, Purdue University, West Lafayette, IN 47907, USA}%
\author{Ivan C.\ Christov}\thanks{Author to whom correspondence should be addressed.}%
\email{christov@purdue.edu}%
\homepage{http://tmnt-lab.org}%
\affiliation{School of Mechanical Engineering, Purdue University, West Lafayette, IN 47907, USA}%

\date{\today}

\begin{abstract}
Low-dimensional dynamical systems are fruitful models for mixing in fluid and granular flows. We study a one-dimensional discontinuous dynamical system (termed ``cutting and shuffling'' of a line segment), and we present a comprehensive computational study of its finite-time mixing properties including the effect of diffusion. To explore a large parameter space, we introduce fit functions for two mixing metrics of choice: the number of cutting interfaces (a standard quantity in dynamical systems theory of interval exchange transformations) and a mixing norm (a more physical measure of mixing). We compute averages of the mixing metrics across different permutations (shuffling protocols), showing that the latter averages are a robust descriptor of mixing for any permutation choice. If the decay of the normalized mixing norm is plotted against the number of map iterations rescaled by the characteristic e-folding time, then universality emerges: mixing norm decay curves across all cutting and shuffling protocols collapse onto a single stretched-exponential profile. Next, we predict this critical number of iterations using the average length of unmixed subsegments of continuous color during cutting and shuffling and a Batchelor-scale-type diffusion argument. This prediction, called a ``stopping time'' for finite Markov chains, compares well with the e-folding time of the stretched-exponential fit. Finally, we quantify the effect of diffusion on cutting and shuffling through a P\'eclet number (a dimensionless inverse diffusivity), showing that the system transitions more sharply from an unmixed initial state to a mixed final state as the P\'eclet number becomes large. Our numerical investigation of cutting and shuffling of a line segment in the presence of diffusion thus present evidence for the latter phenomenon, known as a ``cut-off'' for finite Markov chains, in interval exchange maps.
\end{abstract}

\maketitle

\section{Introduction}
\label{sec:intro}

Even simple discontinuous dynamical systems can exhibit highly nontrivial dynamics and mixing behaviors. One recently studied class of such systems are \emph{piecewise isometries} (PWIs) \cite{g00,g02,d06}. Unlike the \emph{stretching and folding} mechanism of chaotic fluid mixing \cite{o89}, which is underpinned by horseshoe dynamics \cite{o89,sow06} and is sometimes provably the ``best'' mixing possible \cite{wo04}, \textit{cutting and shuffling} underlies granular mixing  \cite{jlosw10,col10,col11,jcol12}. Cutting and shuffling, much like the ``mixing'' of a deck of cards, involves breaking apart the material being mixed into discrete pieces and then putting it back together in a length-preserving ({\it i.e.}, isometric) way \cite{clo11}. Mixing by cutting and shuffling via PWIs on non-Euclidean spaces ({\it e.g.}, the surface of a hemispherical shell) remains an active topic of research \cite{puol16,spuol17}. Meanwhile, fluid mixing by stretching and folding dynamics in physical space is, by now, well-understood \cite{o89,sow06,ko11}. The interaction between stretching and folding and cutting and shuffling, on the other hand, remains a research frontier in the field of dynamical systems. As evidence for the latter point, we refer the reader to the detailed studies by Smith {\it et al.}~\cite{srlm16,srlm17,suol17b} of shear maps coupled to discontinuous motions (such as ``slip deformations''), showing a wealth of distinct types of dynamical behaviors including enhanced mixing and exotic bifurcations.

The simplest example of cutting and shuffling is the one-dimensional PWI known as an \emph{interval exchange transformation} (IET) \cite[\S14.5]{kh95}. Recently, a class of IETs relevant to granular mixing and their mixing properties were studied by Krotter {\it et al.}~\cite{kcol12}. This IET construction and its requisite simulation methodology were introduced to model the intuitive process of cutting and shuffling a line segment. Yu {\it et al.}~\cite{yuol16} extended the work in \cite{kcol12} to account for possible uncertainty in the location of cuts along the line segment, as might be the case when fractionating a granular material such as a powder. Consequently, the length of each portion of the cut and shuffled line segment is random, potentially leading to enhanced mixing. Most recently Smith {\it et al.}~\cite{suol18} introduced a new metric to quantify mixing by IETs, combining the length of the largest uncut subsegment and the evenness of color distribution across subsegments. They demonstrated that cutting  the longest unmixed subsegment of distinct color in half at each iteration, which proves to be computationally inexpensive, can lead to optimal mixing in the sense of minimizing the proposed mixing metric. Such ``optimal'' mixing can also be achieved with a fixed shuffling protocol if the cut locations change at each iteration (beyond just the addition of uncertainty in cut locations as in \cite{yuol16}). Beyond the computational work of \cite{kcol12,yuol16,suol18}, the mathematical theory of IETs is, in fact, quite daunting. The ``weak mixing'' properties of IETs were only recently established \cite{su05,af07} in work that required the development of abstract mathematical notions at such a high level that one of the authors of \cite{af07} (A.~Avila) received the 2014 Fields Medal for his contributions to dynamical systems theory  \cite{fields2014}.

The IET construction introduced by Krotter {\it et al.}~\cite{kcol12} has several parameters that can be varied to produce distinct mixing behaviors, including pathological poor-mixing cases that were examined in detail therein. Some connections between abstract mathematical notions of mixing and numerical experiments with IETs were also summarized in \cite{kcol12,yuol16}, leading to several basic ``design principles'' for how to best cut and shuffle a line segment. Specifically, the number of cuts (subsegments, $N$, introduced in each cutting step) can be varied, the shuffling order (a permutation, $\Pi$, of the integers up to $N$) can be changed, and the lengths of each subsegment (parametrized by a fixed adjacent subsegment length ratio $r$) can be chosen so as to enhance mixing. In \cite{kcol12,yuol16}, mixing was quantified through the percentage of the line segment's length that is unmixed ({\it i.e.}, the proportion of the line segment composed of the same continuous ``color'') or as measured by the number of cutting interfaces ({\it i.e.}, interfaces between different ``colors'' present). Yu~{\it et al.}~\cite{yuol16} additionally showed that even slightly perturbing the cut locations stochastically can break periodicity in the system, again leading to mixing.

Our study of cutting and shuffling is also motivated, in part, by the shuffling of a deck of cards. Numerical results by Trefethen and Trefethen \cite{tt00} illustrate a well-known phenomenon of ``cut-offs'' in card shuffling. Specifically, it is known from work by Diaconis {\it et al.}~\cite{ad86,d95,dsc06}, based on probability theory of finite Markov chains, that it takes about seven riffle shuffles to randomize a fifty-two card deck (see also the {\it New York Times} article \cite{K90} on this fascinating result) \footnote{A simple counting argument by Keller \cite{k95} provides a similar estimate of how the cut-off number of shuffles scales with the number of cards.}. Any further shuffling does not significantly improve the ``mixedness'' of the deck. Hence, seven shuffles represents a \emph{stopping time} for the shuffling process, and the dynamical system exhibits a \emph{cut-off} there, at which a sudden change in the mixedness of the cards from poorly to well mixed occurs; cf.~\cite[Fig.~2]{ad86} and the discussion thereof.

Recently, such cut-offs were demonstrated \cite{lw08} through the numerical simulation of chaotic fluid mixing in a staggered-herringbone microfluidic channel \cite{sdamsw02}. By varying the P\'eclet number (an inverse dimensionless diffusivity), an appropriate global measure of ``mixedness'' can be rescaled and shown to fit the notion of a cut-off, just like in card shuffling. To this end, in this work, we incorporate diffusion into the one-dimensional cutting and shuffling process, and we explore the existence of cut-offs in this model system. Understanding such admittedly ``simple'' systems that yet exhibit ``complex'' dynamical behavior can often be impactful for statistical \cite{krbn10} and material physics \cite{mo14}.

Though there have been a number of theoretical and computational studies of IETs, there is still much to be understood about the basic cutting and shuffling model discussed above. Just as Ashwin {\it et al.}~\cite{ank02} pointed out that ``the mixing properties of interval exchange maps are very subtle and relatively poorly understood and depend on parameters in a sensitive way,'' our aim is to explore how different parameters influence the system's mixing behaviors and whether there is any ``universality'' in the phenomenon.

Going beyond the simple IETs of Ashwin {\it et al.}~\cite{ank02}, Sturman \cite{sturman12} provides a comprehensive review of mathematical theories relevant to discontinuous mixing. Examples of IETs with diffusion were discussed in the context of simulating the mixing of a black-and-white line segment. Mixing with diffusion alone or by cutting and shuffling alone was found to be much less effective than when the two are combined. Going further, Froyland {\it et al.}~\cite{fgw16} proposed local perturbations to speed up mixing in various dynamical systems. Specifically, they sought to optimize how diffusion is added across the system, including one-dimensional line segments undergoing cutting and shuffling. Dynamical systems with different fixed diffusion protocols (termed ``none,'' ``uniform,'' ``Gaussian,'' and ``optimal'') were compared \cite{fgw16}, showing that optimizing the diffusion protocol leads to enhanced mixing at any P\'{e}clet number. Here, we consider only Gaussian diffusion.  

Most recently, Kreczak \textit{et al.}~\cite{ksw17} studied a one-dimensional model of mixing of a line segment with a combination of stretching, permutations and diffusion. Their results show that the global mixing rate depends on both the choice of permutation and the diffusion coefficient. Contrary to expectation (and the results of Ashwin {\it et al.}~\cite{ank02}), increasing the diffusion coefficient leads to a \emph{deceleration} of the mixing rate when \emph{both} stretching and folding and cutting and shuffling are present. Given just four detailed studies \cite{ank02,sturman12,fgw16,ksw17} on this topic exist, the dynamics of cutting and shuffling a line segment in the presence of diffusion remain largely unexplored.

Our work aims to fill a knowledge gap in this field. After introducing the IET construction in Sect.~\ref{sec:IET}, which is the basis of our cutting and shuffling model, we proceed to discuss in detail the effects of diffusion on cutting and shuffling in Sect.~\ref{sec:results}. Although many measures of mixing exist (see, {\it {\it e.g.}}, the detailed review \cite{xlxx17}), including multiscale mixing norms \cite{jlt12}, we use the definition of mixing norm introduced in \cite{ank02} (Sect.~\ref{sec:measure_mix}), which determines how far the line segment's color distribution is from the uniform average color of the initial condition.  Then, having quantified mixing, we verify that diffusion generically leads to decay of the mixing norm  (Sect.~\ref{sec:mixing_diffusion}). To explore the parameter space of this system and see how mixing proceeds under different protocols, we introduce a fit function for the decay of the mixing norm with the number of iterations (Sect.~\ref{sec:measures}). Specifically, from the fit function, we extract a decay time constant for each protocol in the parameter space. Next, we define a time scale, which quantifies the number of iterations for the mixing norm to decay by a factor of $e^{-1}$ (usually termed the e-folding time). With this time scale in hand, in Sect.~\ref{sec:cutoff} we rescale the concentration/color mixing norm decay curves for  different cutting and shuffling systems with diffusion to show that a universal mixing behavior exists. Then, in Sect.~\ref{sec:predict_TPe}, we predict this critical number of iterations using the average length of unmixed subsegments of continuous color during cutting and shuffling and a Batchelor-scale-type diffusion argument. We argue that, on the basis of these numerical results, there is preliminary evidence for the existence of a cut-off phenomenon for IETs with diffusion. Finally, conclusions and avenues for future work are stated in Sect.~\ref{sec:conclusion}. 

\section{\label{sec:IET}SIMULATION METHODOLOGY}

In this section, we describe the simulation methodology that we employ to study mixing of a line segment by cutting and shuffling, including incorporating the effect of diffusion. First, we define a class of \emph{interval exchange transformations} (IETs) that represent our cutting and shuffling protocols.

\subsection{\label{sec:parameters}Parameters of the model}

\begin{figure}[h!]
  \centering
  \includegraphics[width=\columnwidth]{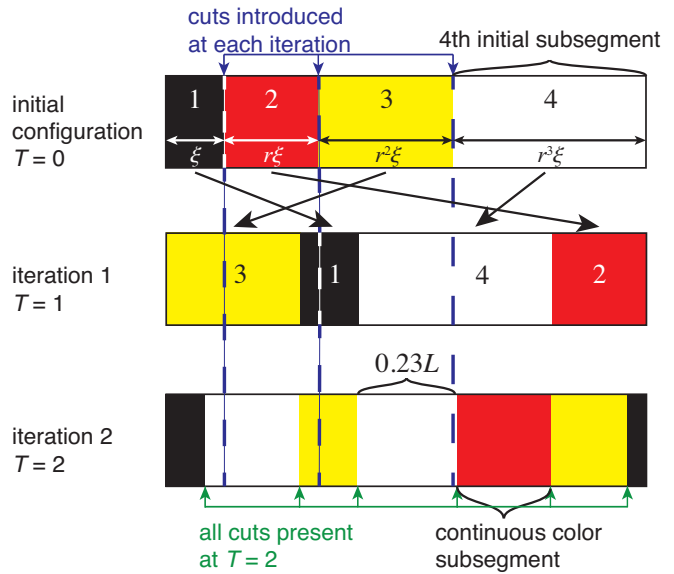}
  \caption{(Color online.) Schematic of how a line segment is cut into $N=4$ (in this case) line segments with adjacent length ratio $r$ and initial subsegment length $\xi$. Two iterations of the cutting and shuffling process (without diffusion) are performed with the permutation $\Pi = [3142]$. Key terminology is labeled. The color map is arbitrary, and the color values are normalized between $0$ and $1$.}
  \label{fig:interval}
\end{figure}

The IET construction can be realized as shown in Fig.~\ref{fig:interval}. The behavior of this dynamical system sensitively depends on three parameters: the number of initial subsegments $N$, the length ratio between adjacent subsegments $r$ and the shuffling permutation $\Pi$. The top row of Fig.~\ref{fig:interval} shows how a line, of total length $L$, is divided into $N=4$ subsegments at each iteration of this dynamical system, which represents the cutting process. Given this value of $N=4$, the line is cut into $N$ pieces and each piece is, additionally, given a distinct color in the initial configuration. The color may, for example, represent different kinds of materials, or the same material but with different ``concentration'' of some tracer being mixed by cutting and shuffling. The ratio $r$ is defined as the ratio of the lengths of adjacent subsegments, while $\xi$ is the length of the first subsegment. Both $r$ and $\xi$ are assumed to be constant in this construction. The permutation $\Pi$ determines the rearrangement order, which represents the shuffling process. Figure~\ref{fig:interval} shows a specific example with $N=4$, $r=1.5$ and $\Pi=[3142]$ \footnote{Here, we use the notation ``$[3142]$'' to denote the permutation that maps the integers $[1234]$ to $[3142]$ in that order. We do not consider $N\ge 10$, so there is no fear of confusion in dropping the spaces between the integers in our notation.}. The cutting location remains the same at every iteration and subsegments are rearranged according to the same pattern. The iteration counter is denoted by $T$, while $T_\mathrm{max}$ is the total number of iterations of cutting and shuffling performed. 

As can be observed from the bottom row in Fig.~\ref{fig:interval}, as the number of iterations $T$ increases, the number of subsegments of continuous color varies, often increasing. Figure~\ref{fig:interval} also highlights the two measures of mixing discussed in \cite{kcol12,yuol16}: the percent unmixed $U(T)$ ({\it i.e.}, the percent out of $L$ corresponding to the longest continuous color subsegment, here $U=23\%$ at $T=2$) and the number of distinct cuts $C(T)$ ({\it i.e.}, interfaces between colors, here $C = 6$ at $T=2$).

\subsection{\label{sec:choose_params}Choosing the parameters: Design rules}

A major conclusion of Krotter {\it et al.}~\cite{kcol12} was that mixing under this type of IET reaches a ``point of diminishing returns'' as $N$ increases with four to five subsegments being sufficient to produce significant shuffling of the material. Thus, in the present work, we restrict our attention to the cases $N=4$ and $N=5$.

As discussed in \cite{kcol12} and further elucidated in \cite{yuol16}, in choosing the shuffling permutation $\Pi$, we should exclude ones that lead to pathological behaviors. Specifically, we only consider permutations $\Pi$ that are (i) irreducible, (ii) non-rotational, (iii) without the first or last element fixed, and \emph{additionally} for the cases with $N>3$, (iv) without a number of elements $> 1$ and $\le N-2$ being consecutive. For example, if $N=4$, there are nine ``allowed'' permutations: $\Pi = [2 4 1 3]$, $[2 4 3 1]$, $[3 1 4 2]$, $[3 2 4 1]$, $[3 4 2 1]$, $[4 1 3 2]$, $[4 2 1 3]$, $[4 3 1 2]$, $[4 3 2 1]$. In particular, $\Pi = [2 1 4 3]$ is excluded by rule (i) for being reducible, {\it i.e.}, elements $1$ and $2$ are interchanged and $4$ and $3$ interchanged, splitting (``reducing'') the permutation into two sub-permutations. Meanwhile, $\Pi = [2 3 4 1]$ is excluded by rule (ii) for being a rotation, {\it i.e.}, elements are shifted one to the right without significant re-arrangement. Then, $\Pi = [4 2 3 1]$ is excluded by rule (iv) because elements and $2$ and $3$ remain consecutive in the permutation.

Additionally, consonant with the available theory of IETs \cite{k75}, Krotter {\it et al.}~\cite{kcol12} showed that the initial cuts should break apart the interval in such a way that the adjacent segment length ratio $r$ is ``closer'' to an irrational number. Moreover, it was concluded in \cite{kcol12} that the initial distribution of subsegment lengths should be ``balanced,'' that is, $r$ should be chosen close (but not equal) to unity. Thus, in the present work, without loss of generality, the ratio is taken to be $r>1$, then as shown in Fig.~\ref{fig:interval}, the total length $L$ of the line segment is
\begin{equation}
L = \sum _{j=1}^{N} r^{j-1} \xi.
\label{eq:L_x}
\end{equation}

In order to realize the shuffling process within a MATLAB code, while ensuring length-preservation of the line segment without being subject to round-off errors, the initial lengths of \emph{all} subsegments, {\it i.e.}, $\{\xi, r\xi, r^2\xi, \hdots, r^{N-1}\xi\}$, should \emph{all} be integers. If the latter condition is met, then cuts \emph{always} fall at unique indexes in the discrete array that represents the line segment computationally, which ensures that no length can be ``lost'' by rounding potentially fractional indexes up or down. To restate this important point: by guaranteeing that cuts fall at unique array elements, then we automatically ensure that the line segment's length cannot change, {\it i.e.}, we enforce conservation of mass. 

To ensure that the subsegment lengths $\{\xi, r\xi, r^2\xi, \hdots, r^{N-1}\xi\}$ are all integers, we convert $r$ to a fraction as $r_n/r_d$, where $r_n$ is the nominator (an integer), and $r_d$ is the denominator (also an integer). This conversion is always possible as long as $r$ is a rational number. Since, in MATLAB, we use finite-precision floating point arithmetic, then any $r$ we could pick must be representable as a fraction, though it might be quite a ``complicated'' fraction. Now, once we have written $r=r_n/r_d$, it is clear that multiplying the list $\{1, r, r^2, \hdots, r^{N-1}\}$ by $r_d^{N-1}$ yields a set of integers with greatest common divisor $1$. Thus, we conclude that
\begin{equation}
\xi = r_{d}^{N-1} 
\label{eq:x_int}
\end{equation}
will guarantee that every subsegment's length is an integer. Finally, substituting the expression for $\xi$ from Eq.~\eqref{eq:x_int} into the total length of the line given in Eq.~\eqref{eq:L_x}, we find that
\begin{equation}
L = \sum _{j=1}^{N} r^{j-1} r_{d}^{N-1} = \left( \frac{1-r^N}{1-r} \right) r_d^{N-1}.
\label{eq:L}
\end{equation}
In particular, we conclude that $L$ is an integer, and it depends only on $N$ and $r$ for a given cutting and shuffling protocol of our integer-arithmetic IET construction. In plots, we will generally normalize the horizontal axis by $L$ so that the line segment $[0,L]$ becomes the interval $[0,1]$, and subsegments' positions along this unit interval are displayed.

\subsection{\label{sec:measure_mix}Quantifying mixing: Cutting interfaces and the mixing norm}
After $T$ iterations of the cutting and shuffling map, the initially coherent sets of colors assigned to the pieces of the line segment can form a complex and intricate pattern \cite{kcol12,yuol16}. There are many measures of mixing that can be used to quantify the \emph{degree of mixing} produced by cutting and shuffling. On the one hand, there are discrete measures such as counting the number of distinct cuts between different colors, as discussed earlier. On the other hand, there are a variety of so-called ``mixing norms'' that can be used to quantify the degree of mixing in a more ``global'' way, as discussed in the review by Thiffeault \cite{jlt12}. In this work, we will utilize the number of cuts $C(T)$ and the $L^p$ function-space norm, denoted $||c||_{p}(T)$, to quantify mixing. 

As discussed in Sect.~\ref{sec:parameters}, the number of cutting interfaces $C(T)$ refers to the number of distinct interfaces between two different continuous color subsegments after $T$ iterations. This is a clear and intuitive metric of how much rearrangement our cutting and shuffling protocol has achieved. Taking Fig.~\ref{fig:interval} as an example, after the first iteration, there are 3 interfaces since 4 subsegments are generated, therefore $C(1)=3$. After the second iteration, there are 6 interfaces between distinct components, therefore $C(2)=6$. Given a color distribution $\{c_i\}_{i=1}^L$ across the line segment (represented by the lattice $i=1,2,\hdots,L$) after $T$ iterations, we can explicitly define the number of cutting interfaces as
\begin{equation}
C(T) = \sum_{i=1}^{L-1} \lceil c_{i+1} - c_i \rceil,
\label{eq:new_cut_def}
\end{equation}
since $0\le c_i \le 1$ by construction. Thus, the number of cuts is not a true proxy of mixing as it does not take into account whether or how the color changes on average. Nevertheless, the number of cuts is meaningful because, if reactions and diffusion are added into the cutting and shuffling protocol, then it is expected that having a wide distribution of cuts will lead to fast homogenization of the material \cite{mo89,ccr99}. The number of cuts (discontinuities) is also a quantity of interest in the abstract mathematical discussion of IETs \cite{n09}. 

To mitigate some of the weakness of $C(T)$ as a measure of mixing, we also use a mixing norm to quantify mixing. Specifically, we define an ``$L^p$ norm'' of a function $c(X,T)$ as
\begin{equation}
||c||_{p}(T) = \left(\frac{\int_0^L |c(X,T) - \bar{c}|^p \,\mathrm{d}X}{\int_0^L \mathrm{d}X}\right)^{1/p},
\end{equation}
where $1\le p < \infty$ and
\begin{equation}
\bar{c} = \frac{\int_0^L c(X,T) \,\mathrm{d}X}{\int_0^L \mathrm{d}X}
\label{eq:mix_norm_define}
\end{equation}
is the \emph{average color} of the line segment. Here, $X$ is continuous variable running along the length of the line segment: $0\le X \le L$. For concreteness, when initially constructing the line segment, we assign each subinterval $i$ a color value $(i-1)/(N-1)$ ($i=1,\hdots,N$), which is always between $0$ and $1$. In the plots presented below, we can use various color maps in MATLAB to make the colors stand out visually. Thus, the mixing norm $||c||_{p}(T)$ measures how far the segment's color distribution is from the uniform average concentration/color $\bar{c}$, in an appropriately global way. The case of $p=2$ is of interest as it measures the variance, or root-mean-square deviation, of the color distribution. As discussed by Thiffeault \cite[p.~5]{t08}, ``[v]ariance is thus a useful measure of mixing: the smaller the variance, the better the mixing.'' And, while the average $\bar{c}$ remains constant in time and unchanged by diffusion for periodic boundary conditions, the variance is depleted as long as there is diffusion and non-zero color gradients along the line segment \cite{t08}, with the cutting and shuffling process controlling how gradients are created and, thus, the decay rate.

\begin{figure*}[ht]
  \centering
  \includegraphics[width=\textwidth]{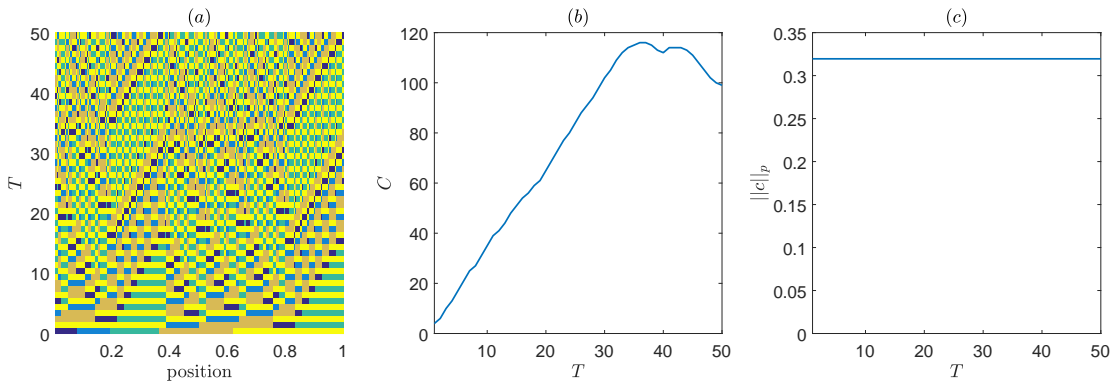}
  \caption{(Color online.) Mixing by cutting and shuffling with $N=5$, $r=1.5$, $\Pi = [5 2 4 1 3]$, $T_\mathrm{max}=50$. (a) Space-time plot of the color evolution. (b) Growth of the number of cutting interface $C(T)$. (c) Evolution of the mixing norm $\|c\|_p(T)$.}
  \label{fig:space-time} 
\end{figure*}

To compute the mixing norm from the discrete data of our cutting and shuffling simulations, consider a distribution of colors $c(X,T)$ across the line segment $0\le X\le L$. After $T$ iterations of cutting and shuffling, $c_{j}$ represents the color value of the $j$th continuous-color piece, where $j=1,\hdots,k(T)$, and $k(T)$ is the number of continuous-color subsegments after $T$ iterations. Then, we can compute the integrals in the definition of the mixing norm exactly over each continuous-color piece and reduce the definition from Eq.~\eqref{eq:mix_norm_define} to
\begin{equation} 
	||c||_{p}(T) = \left(\frac{ \sum_{j=1}^{k(T)} |c_{j}-\overline{c}|^{p} l_{j}}{\sum_{j=1}^{k(T)} l_j} \right)^{{1}/{p}},
\label{eq:mix_norm}
\end{equation}
where
\begin{equation}
\bar{c} = \bar{c}(r,N) = \frac{\sum_{j=1}^{N} c_{j} l_{j}}{\sum_{j=1}^{N} l_{j}}  \qquad (T = 0)
\label{eq:cbar}
\end{equation}
is, as before, the uniform average concentration/color of the initial condition. Note that $\sum_{j=1}^{k(T)} l_j=L$ by definition, where $L$ is given by Eq.~\eqref{eq:L}, and $l_{j}$ is the length of $j$th subsegment of continuous color. At iteration $T$, there are $1\le k(T) \le L$ pieces of continuous color with  $k(0)=N$. For $p=2$, Eq.~\eqref{eq:mix_norm} describes the standard deviation (square root of the variance) of the mixture's ``concentration.'' If normalized appropriately, the $p=2$ mix norm can be made to agree with Dankwerts' classical definition of the \emph{intensity of segregation} \cite{d52}, which he used to quantify mixing. Henceforth, we restrict to the $p=2$ case for the remainder of this work.

\subsection{\label{sec:st_plots}Visualizing mixing: Space-time plots}

Mixing of a line segment by cutting and shuffling can be visually represented by space-time plots as initially discussed in \cite{kcol12}. To create a space-time plot, we combine all line segments after each iteration and arrange them from bottom to top in a two-dimensional (2D) space-time plot. The horizontal axis is the dimensionless lattice position $X/L$, while the vertical axis represents the number of iterations $T$. From space-time plots, we can visually identify the evolution of mixing, including periodic behavior and poor mixing. Figure~\ref{fig:space-time}(a) shows an example of a space-time plot of a cutting and shuffling protocol. 

In Fig.~\ref{fig:space-time}(b), we plot the number of cuts $C(T)$, which grows in time but levels out after a while and starts decreasing. This protocol is ultimately periodic, so even over many iterations it does not produce good mixing. In Fig.~\ref{fig:space-time}(c), we plot the mixing norm $\|c\|_p(T)$ for this case. The mixing norm remains constant, meaning that the color distribution never approaches the average color. Of course, since cutting and shuffling merely redistributes the color pieces, without changing their individual colors, the distribution cannot approach the average (see also \cite[Sect.~4]{fgw16}). The latter is, of course, the classical distinction of \emph{stirring} versus \emph{mixing} \cite{o89sa,o89,w88}. That is to say, while cutting and shuffling (the ``mechanical'' stirring process in our approach) can significantly disperse the initially continuous color segments, diffusion (such as molecular diffusion in a fluid or collisional diffusion in a granular flow) is needed to ultimately homogenize and mix the material. This distinction and the interplay between stirring and mixing bring us to a key contribution of the present work: incorporating diffusion into IETs and examining the resulting universal mixing behaviors.

\subsection{\label{sec:diffusion}Incorporating diffusion}

To ensure complete and thorough mixing of a line segment, {\it i.e.}, $c(X,T) \to \bar{c}$ for all $X\in[0,L]$ as $T\to\infty$, we must incorporate diffusion into the system. As is well known, diffusion by itself would mix an initial line segment (such as the one shown in the top row of Fig.~\ref{fig:interval}) very slowly. Thus, here we are interested in the \emph{nontrivial} interaction of cutting and shuffling (redistribution of color) and diffusion (relaxation of the color distribution to the mean).

Following Pierrehumbert \cite{p94}, we would like to incorporate a \emph{time-discrete} diffusion step between cutting and shuffling steps. To this end, consider a generic diffusion equation for the concentration/color $c(X,T)$ with characteristic diffusivity $D$:
\begin{equation}
\frac{\partial c}{\partial T} = D\frac{\partial^2 c}{\partial X^2},
\label{eq:diffusion_eq}
\end{equation}
where $X\in[0,L]$ and $T\in[0,T_\mathrm{max}]$. We can discretize Eq.~\eqref{eq:diffusion_eq} using the usual forward-time, central-space (FTCS) scheme:
\begin{equation}
\frac{c_i^{n+1}-c_i^n}{\Delta T} = D\frac{c^n_{i+1}-2c^n_i+c^n_{i-1}}{(\Delta X)^2},
\label{eq:diffusion_scheme}
\end{equation}
where $c^n_i\approx c(X_i,T^n)$ with $X_i = i\Delta X$ and $T^n = n\Delta T$. This discretization is stable if $D\Delta T \le \tfrac{1}{2}(\Delta X)^2$ \cite[\S6.3]{Strikwerda}.

However, in our cutting and shuffling protocols, the color is only defined on integer lattice points, and we iterate by integer increments in time, {\it i.e.}, $\Delta X = \Delta T = 1$. Then, given the standard stability criterion for the finite-difference scheme in Eq.~\eqref{eq:diffusion_scheme}, we must restrict our equivalent diffusion coefficient such that $D \le {1}/{2}$.

Therefore, diffusion can be incorporated into cutting and shuffling by an additional sweep through the lattice at each iteration of the IET. That is to say, after the line segment is cut and shuffled as described above, an additional sweep through the lattice points is performed using the replacement rule:
\begin{equation} 
	c_{i} \mapsto (1-2D)c_{i}+Dc_{i+1}+Dc_{i-1},
\label{eq:diffusion_rule}    
\end{equation}
where $c_{i}$ is the color value at the $i$th lattice point ($i=1,2,\hdots, L$). A common choice is $D=1/2$, in which case our replacement rule from Eq.~\eqref{eq:diffusion_rule} becomes a simple averaging: $c_{i} \mapsto \frac{1}{2}\left( c_{i+1} + c_{i-1}\right)$. Given a particular diffusion coefficient $D$, we would like to show that the line segment can be completely mixed in many fewer iterations than by the IET alone. To completely specify the diffusion rule, periodic boundary conditions are used to set $c_{L+1} = c_1$ and $c_0 = c_L$.

\begin{table}
\centering
\begin{tabular}{l l l l l}
$r$ & $r_n$ & $\xi$ & $L$ & $T_\mathrm{max}$ \\
\hline
1.25 & 5 & 64 & 369 & 50 \\
1.2 & 6 & 125 & 671 & 166 \\
1.4 & 7 & 125 & 888 & 290 \\
1.6 & 8 & 125 & 1,157 & 492 \\
1.8 & 9 & 125 & 1,484 & 809 \\
1.1 & 11 & 1,000 & 4,641 & 7,910 \\
1.3 & 13 & 1,000 & 6,187 & 14,057 \\
\hline
\end{tabular}
\caption{Dependence of $L$ on $r$, here specifically for $N=4$, and the corresponding maximum number of iterations, via Eq.~\eqref{eq:3}, with the first row ($r=1.25$ and $T_\mathrm{max}=50$) being the reference system for the remaining.}
\label{tb:rLT}
\end{table}

We would like all of our simulations to have the same ``effective'' diffusion coefficient, $D$, as it might arise from inter-particle collisions in a granular medium (see, {\it e.g.}, \cite{fuol15}). However, the length $L$ of the lattice depends on $N$ and $r$ as discussed in Sect.~\ref{sec:choose_params}, see Table~\ref{tb:rLT}. Hence, a fixed diffusion coefficient does not yield the same behavior on different lattices, over the same number of iterations, simply because of the change in the line segment's length. To ensure an ``equivalent'' diffusion behavior across lattices of different $L$, we apply dimensional analysis to connect the number of iterations $T_\mathrm{max}$ that a given IET is required to run for, given a fixed diffusivity $D$ but different domain lengths. The problem reduces to matching dimensionless diffusion coefficients once the domain is mapped from $[0,L]$ to $[0,1]$ and the number of iterations (from $0$ to $T$) is normalized to a discrete time-like variable running from $0$ to $1$. 

To this end, consider two cases of cutting and shuffling protocols with diffusion, the first with $D_1$, $T_{\mathrm{max},1}$ and $L_1$ and the second with $D_2$, $T_{\mathrm{max},2}$ and $L_2$. From dimensional analysis, we should ensure that the dimensionless diffusion coefficients ({\it i.e.}, the inverse P\'eclet numbers assuming the intrinsic ``velocity'' scale $L/T$) match:
\begin{equation}
\frac{D_2 T_{\mathrm{max},2}}{L_2^2}= \frac{D_1 T_{\mathrm{max},1}}{L_{1}^{2}}.
\label{eq:dimensional analysis}
\end{equation}
Assuming equal diffusivity ($D_1=D_2$), the diffusion coefficient cancels out, and we can relate the number of iterations $T_{\mathrm{max},2}$ needed on a lattice of length $L_2$ to those ($T_{\mathrm{max},1}$ and $L_1$) of the reference lattice:
\begin{equation}
T_{\mathrm{max},2} = \left(\frac{L_{2}}{L_{1}}\right)^{2} T_{\mathrm{max},1}. 
\label{eq:3}
\end{equation}

\section{\label{sec:results}RESULTS AND DISCUSSION}

In this section, we examine the mixing outcomes of cutting and shuffling in the presence of diffusion. Specifically, we address the hypothesis that cut-offs and universal behavior exist in the family of IETs with diffusion that we have described/constructed above. 

\begin{figure*}[ht]
  \centering
  \includegraphics[width=\textwidth]{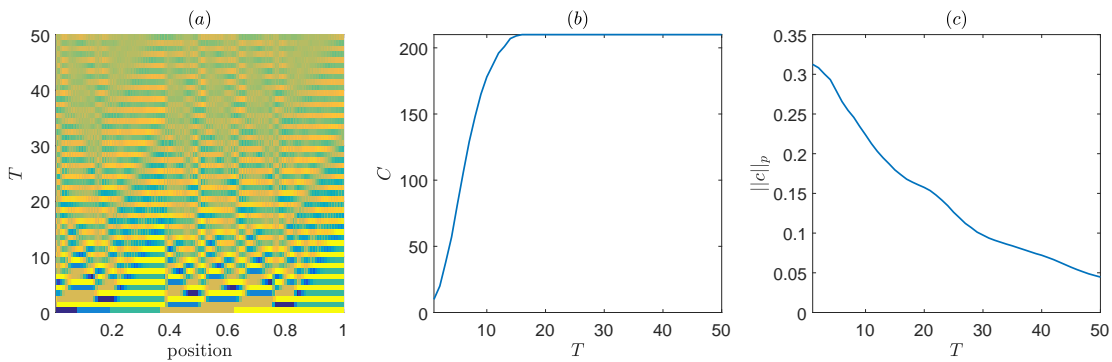}
  \caption{(Color online.) Incorporating diffusion with $D=0.5$ into the IET from Fig.~\ref{fig:space-time} with $N=5$, $r=1.5$, $\Pi = [5 2 4 1 3]$, and $T_\mathrm{max}=50$.}
  \label{fig:diffusion-space-time}
\end{figure*}

\begin{figure*}[ht]
  \centering
  \includegraphics[width=\textwidth]{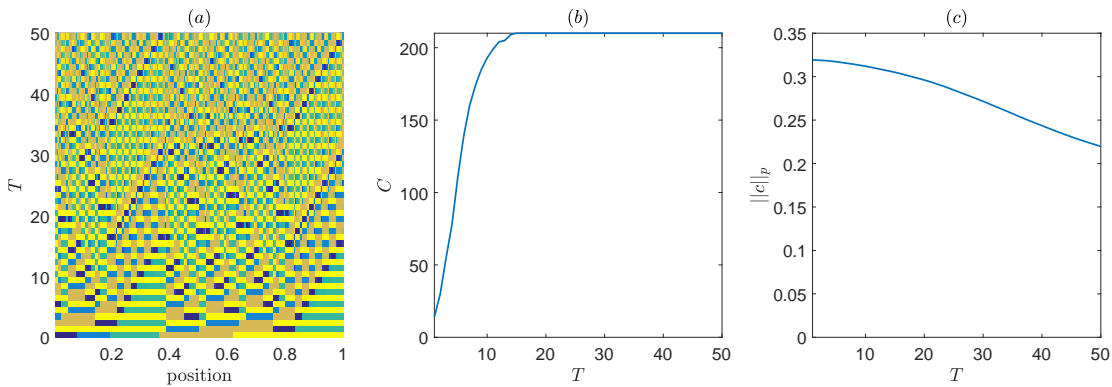}
  \caption{(Color online.) Incorporating diffusion with $D=0.01$ into the IET from Fig.~\ref{fig:space-time} with $N=5$, $r=1.5$, $\Pi = [5 2 4 1 3]$, and $T_\mathrm{max}=50$.}
  \label{fig:small-diffusion-space-time}
\end{figure*}

Previous work has sufficiently addressed the non-diffusive ({\it i.e.}, deterministic) mixing by IETs both mathematically and through simulation studies. Therefore, in this subsection, we summarize just the key results. From \cite{kcol12,yuol16}, it is clear that the number of subsegments $N$, the permutation $\Pi$ and the ratio $r$ greatly affect the mixing outcomes. Specifically, three ``design rules'' have been suggested to improve the mixing behavior (recall the discussion in Sect.~\ref{sec:choose_params}):
\begin{enumerate}
\item[(i)] Reducible and rotational permutations, as well as those that exhibit unsatisfactory shuffling, should be excluded.
\item[(ii)] The line segment should be cut into no more than six subsegments ({\it i.e.}, $N\leq6$), larger $N$ do not significantly improve mixing.
\item[(iii)] The ratio $r$ should not be ``large.'' Specifically, it should ideally be an irrational number close to 1. (Continued fraction expansions can be used to quantify ``how irrational'' $r$ is.)
\end{enumerate} 

\subsection{\label{sec:mixing_diffusion}Mixing behavior of IETs with diffusion}

In the present work, motivated by the hypothesis that IETs with diffusion exhibit cut-offs in the sense of card shuffling, we are first interested in establishing how the mixing behaviors previously studied are affected by the presence of diffusion. Specifically, we study the effect of the diffusion coefficient's magnitude ({\it e.g.}, $D = 0, 0.01, 0.5, \hdots$), having matched the total iterations $T$ using the dimensional analysis rules from Section~\ref{sec:diffusion} to ensure comparable ``amounts'' of diffusion across lattices of different lengths.

As an introductory example, let us consider how a typical space-time plot, such as the one shown in Fig.~\ref{fig:space-time}(a) changes when diffusivity with $D=0.5$ is incorporated into the cutting and shuffling process. Under the same parameters as Fig.~\ref{fig:space-time}(a), Fig.~\ref{fig:diffusion-space-time}(a) shows the space-time plot of mixing by cutting and shuffling including diffusion. The most obvious effect is that the space-time plot becomes ``fuzzy'' as diffusion now blurs the different colors of the subsegments being cut and shuffled about.

We can also examine how the number of cutting interfaces $C(T)$ and the mixing norm $\|c\|_p(T)$ evolve in the presence of diffusion. Compared with Fig.~\ref{fig:space-time}(b), the number of cutting interfaces in Fig.~\ref{fig:diffusion-space-time}(b) grows quickly and reaches an absolute maximum. This difference in how $C(T)$ evolves is due to the fact that material is now not just cut and shuffled but also mixed by diffusion. Diffusion changes the color of nearby lattice points through the diffusion rule [given by Eq.~\eqref{eq:diffusion_rule}], thereby quickly causing nearby lattice points to have slightly different color values, and all of these slight changes are counted as cutting interfaces by Eq.~\eqref{eq:new_cut_def}. In a small number of iterations, the number of cutting interfaces $C(T)$ increases without exhibiting periodic patterns, and reaches its absolute maximum value $\max_T C(T) = L-1$. This upper bound is due to the fact that eventually the color of \emph{every} lattice point is distinct from every other (even if just slightly so) due to cutting, shuffling and diffusion. 

The mixing norm $||c||_{p}(T)$, on the other hand, now decreases (asymptotically to 0) with $T$, as seen in Fig.~\ref{fig:diffusion-space-time}(c), instead of remaining constant as in Fig.~\ref{fig:space-time}(c). In the presence of diffusion, the cutting and shuffling process eventually drives the color of the line segment to the average one, $\bar{c}$, which is set by the initial conditions.

Next, we would like to establish the effect of varying the diffusion coefficient $D$ ({\it i.e.}, ``small'' diffusivity versus ``large'' diffusivity). For a smaller diffusion coefficient of $D=0.01$, the mixing behavior is shown in Fig.~\ref{fig:small-diffusion-space-time}. The growth of cutting interfaces $C(T)$ in Fig.~\ref{fig:diffusion-space-time}(b) is almost the same as in Fig.~\ref{fig:small-diffusion-space-time}(b), showing a weak sensitivity to the diffusivity. This observation suggests that $C(T)$ might not be an effective way to measure the degree of mixing among systems with different diffusion coefficients. While in the deterministic (no diffusion) case, insightful mathematical results can be obtained about the growth of the number of cuts \cite{n09}, any amount of diffusion perturbs the color values so that differences in color that are counted as ``cuts'' appear immediately. Although one can invent threshold criteria for how much change $|c_{i+1}-c_i|$ should signal a ``cut'' in Eq.~\eqref{eq:new_cut_def}, this is ultimately a fruitless task. Meanwhile the mixing norm $\|c\|_p(T)$ in Fig.~\ref{fig:small-diffusion-space-time}(c) decays much more slowly that in Fig.~\ref{fig:diffusion-space-time}(c), showing (as is to be expected on the basis of previous studies, {\it e.g.}, \cite{ank02,sturman12,fgw16,ksw17}) that the mixing norm effectively quantifies the differences in  mixing processes with ``small'' diffusivity ($D=0.01$) versus ``large'' diffusivity ($D=0.5$).

\subsection{\label{sec:measures} Quantifying the effect of diffusion on the decay of the mixing norm}

Based on the numerical evidence in Figs.~\ref{fig:diffusion-space-time}(c) and \ref{fig:small-diffusion-space-time}(c), we suggest that the evolution of $||c||_{p}(T)$ can be approximated by a \emph{stretched-exponential} function \footnote{The choice of a stretched-exponential function is based on the fact that such a parametrization arises and accurately describes a wide range of relaxation processes in disordered condensed matter systems \cite{p96}.} parametrization:
\begin{equation}
||c||_{p}(T) \simeq M\cdot e^{-\left(T/\tau\right)^\alpha},
\label{eq:cnorm_fit}
\end{equation}
where $M = ||c||_{p}(0)$ is the initial norm before mixing but $\tau$ and $\alpha$ are \emph{a priori} unknown fitting parameters. The time constant $\tau$ quantifies how fast the mixing norm decays with $T$, while $\alpha$ determines how skewed the decay curve is. If $\alpha=1$, the fit function in Eq.~\eqref{eq:cnorm_fit} is a ``perfect'' exponential, while for $\alpha < 1$, it is skewed and decays more slowly (asymptotically as $T\to\infty$). 

\begin{figure}[ht]
  \centering
  \includegraphics[width=0.5\textwidth]{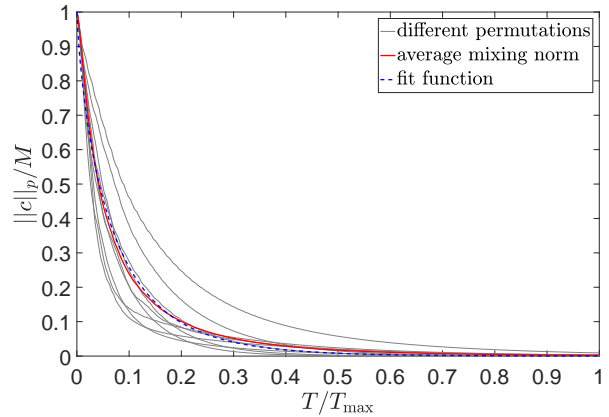}
  \caption{(Color online.) Mixing norm decay curves $||c||_{p}(T)$ with $N = 4$, $r = 1.25$, $M = 0.3650$, $D = 0.5$ and $T_\mathrm{max} = 1000$ iterations for various $\Pi$ (as chosen according to the rules in Sect.~\ref{sec:choose_params}), shown as gray curves. The average mixing norm decay curve is the bold curve (red online), and its fit is the dashed curve (blue online). The best-fit parameters, according to Eq.~\eqref{eq:cnorm_fit} for the average curve, are $\tau = 68.17$ and $\alpha = 0.7866$.}
  \label{fig:fit-norm}
\end{figure}

\begin{figure}[ht]
  \centering
  \includegraphics[width=0.5\textwidth]{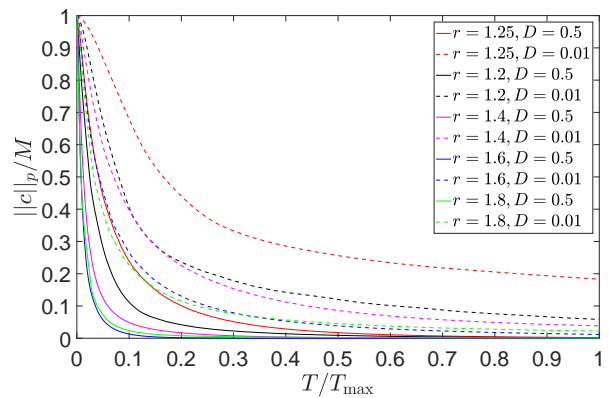}
  \caption{(Color online.) The average (across permutations) mixing norm $||c||_{p}(T)$ decay curve from Fig.~\ref{fig:fit-norm} as $D$ (dashed for $D=0.01$ and solid for $D=0.5$) and $r$ (different colors) are varied. All curves are rescaled with respect to  $T_\mathrm{max}=1000$.}
  \label{fig:different-average-norm}
\end{figure}

For a fixed ratio $r$, number of subsegments $N$ and diffusivity $D$, we average the mixing norm evolution curves $||c||_{p}(T)$ across different permutations, then we fit the averaged profile to Eq.~\eqref{eq:cnorm_fit}. MATLAB's nonlinear least-squares subroutine {\tt lsqcurvefit} is used to obtain the best-fit values of the parameters. Thus, we obtain numerical values for $\tau$ and $\alpha$. Figure \ref{fig:fit-norm} illustrates this procedure for $N = 4$, $r = 1.25$ and $D = 0.5$: the light gray curves represent $||c||_{p}(T)$ for the different permutations $\Pi$ considered, while the dark gray (red online) curve is the average value of $||c||_{p}(T)$ across permutations, and the dashed curve (blue online) is the fit according to Eq.~\eqref{eq:cnorm_fit}. We observe that the proposed fit function captures the overall trend of the decay of the average mixing norm quite well.

Next, we vary the diffusion coefficient $D$ and subsegment length ratio $r$, plotting the fit curves of the average in Fig.~\ref{fig:different-average-norm}. The more ``complex'' $r$  is (see, {\it e.g.}, the discussion in \cite{kcol12} about defining $r$ as a continued-fraction expansion of increasing length), the fewer iterations it takes to reach complete mixing. Generally, though there are exceptions, as a careful examination of Fig.~\ref{fig:different-average-norm} reveals. Nevertheless, it takes fewer iterations to homogenize the line segment with a larger diffusion coefficient ($D=0.5$) than with a smaller diffusion coefficient ($D=0.01$), as expected. As discussed above, the connection between the value of $r$ and the resulting cutting and shuffling behavior is highly nontrivial. Thus, our approach of computing the average mixing curve (over all permutations considered) and fitting it via two parameters provides a clear and quantitative way to compare protocols with different $r$ and $D$. 

Table \ref{table:norm-fit-parameter} summarizes the values of $\tau$ and $\alpha$ for the different cutting and shuffling systems shown in Figs.~\ref{fig:fit-norm} and \ref{fig:different-average-norm}. Each  stretched-exponential fit quite accurately describes an individual averaged mixing curve. Specifically, the values for $\tau$ and $\alpha$ reported to four significant digits in Table \ref{table:norm-fit-parameter} are within the fits' 95\% confidence intervals to two significant digits. Now, from the values of the fitting parameters in Table \ref{table:norm-fit-parameter}, we can infer that the time constant $\tau$ is smaller for systems that reach complete mixing in fewer iterations, as expected. Thus, amongst the five choices of $r$ considered, 
\begin{equation}
r = 1.6 = 1 + \frac{1}{1 + \frac{1}{1 + \frac{1}{2}}} \equiv [1; 1, 1, 2]
\label{eq:r16fract}
\end{equation}
results in the swiftest mixing (smallest time constant $\tau$); the continued fraction expansions for the remaining $r$ values are provided in Table~\ref{table:norm-fit-parameter}.

\begin{table}
\begin{tabular}{l c c c c}
& \multicolumn{2}{c}{$D=0.5$} & \multicolumn{2}{c}{$D=0.01$} \\
& \multicolumn{2}{c}{\hrulefill} & \multicolumn{2}{c}{\hrulefill}\\
$r$ & $\tau$ & $\alpha$ & $\tau$ & $\alpha$ \\
\hline
$1.25 = [1; 4]$ & 68.17 & 0.7866 & 333.4 & 0.6471\\
$1.2 = [1; 5]$ & 39.73 & 0.8118 & 14.11 & 0.6385\\
$1.4 = [1; 2, 2]$ & 20.26 & 0.7772 & 12.26 & 0.6660\\
$1.6 = [1; 1, 1, 2]$ & 12.77 & 0.8526 & 73.15 & 0.6794\\
$1.8 = [1; 1, 4]$ & 13.73 & 0.7708 & 60.43 & 0.6017\\
\hline
\end{tabular}
\caption{Fit parameters of $||c||_{p}(T)$ according to the model in Eq.~\eqref{eq:cnorm_fit} for the \emph{average} mixing norm decay curves shown in Fig.~\ref{fig:different-average-norm}. Continued fraction expansions of $r$ are given using the notation in Eq.~\eqref{eq:r16fract}.}
\label{table:norm-fit-parameter}
\end{table}

\begin{figure}[t]
  \centering
  \includegraphics[width=0.5\textwidth]{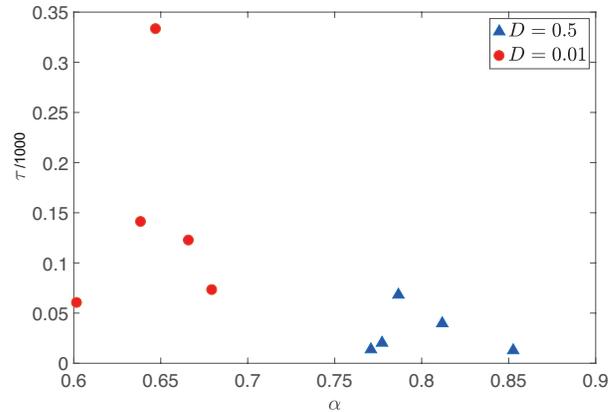}
  \caption{(Color online.) Scatter plot in the $(\alpha,\tau)$ plane of the fit parameters for the average mixing norm decay curves from Fig.~\ref{fig:different-average-norm} (see also Table~\ref{table:norm-fit-parameter}). ``Good'' mixing is observed near the bottom-right corner of the scatter plot, while ``poor'' mixing is observed near the top-left corner.}
  \label{fig:scatter-plot}
\end{figure}

Figure \ref{fig:scatter-plot} shows the scatter plot of $\tau$ and $\alpha$ values. This figure and approach to analyzing our data is inspired by the so-called ``$\tau$-$bias$'' scatter plots of McIlhany and Wiggins \cite{mcwig}. In \cite{mcwig}, the normalized variance of concentration was used to quantify the degree of fluid mixing in a microfluidic device. The parameters $\tau$ and ${bias}$ were introduced to quantify the evolution of the concentration variance curve and, thus, mixing. The parameter $\tau$ is, just as in the present work, interpreted as a time constant describing the temporal decay of the variance of concentration, while $bias$ quantifies the ``unmixedness'' of the final asymptotic state. In our work, $bias=0$ in all cases because the cutting and shuffling process with diffusion leads to $c_i \to \bar{c}$ for every lattice site $i$ as $T\to\infty$. In \cite{mcwig}, it was suggested that small values of \emph{both} $\tau$ and \textit{bias} correspond to ``good'' mixing cases. For our problem, the best mixing case is on the bottom right of Fig.~\ref{fig:scatter-plot}, which corresponds to small $\tau$ but $\alpha$ closer to 1 ({\it i.e.}, not small). The worst mixing cases are on the top left  of the figure, which corresponds to a large $\tau$ and $\alpha$ far away from 1. Therefore, we find out the relationship between fitting parameters and mixing behaviors. Evidently, $\tau$ and $\alpha$ exhibit nontrivial dependences on both the length ratio $r$ and the diffusion coefficient $D$. 

\subsection{\label{sec:cutoff}Universal mixing curves and cut-offs}

So far, we have shown (a) how to incorporate diffusion into the mixing of a line segment by cutting and shuffling and (b) how to quantify mixing across families of protocols via the number of cutting interfaces and the mixing norm. In this section, we would like to substantiate, through numerical results, the central hypothesis of this work: namely that ``cut-offs'' (and the associated concept of ``stopping times'') exist in IETs with diffusion and, furthermore, all mixing behaviors exhibited by cutting and shuffling with diffusion are, in a sense, \emph{universal}. Guided by the work in \cite{tt00,d95,lw08}, we now turn our attention collapsing the mixing norm decay curves of different IETs with diffusion onto a universal profile. As Liang and West \cite{lw08} note, ``[t]o prove the existence of a cutoff is in general very hard, relying on special features of the sequence of systems,'' thus, for the present purposes, we also settle for numerical evidence thereof.

Liang and West \cite{lw08} used the number of iterations, denoted by $T_{Pe}$, required to decrease the initial value of the mixing norm by $50\%$, {\it i.e.}, $T_{Pe}$ such that $||c||_{p}(T_{Pe}) \approx 0.5||c||_{p}(0)$, to collapse the mixing norm curves across different model parameters. While we can certainly compute such a $T_{Pe}$ value numerically from the decay curves of $||c||_{p}(T)$, here it is natural to use the \emph{e-folding time} of the stretched exponential fit from Eq.~\eqref{eq:cnorm_fit}. In other words, we define $T_{Pe}$ as the number of iterations required for $||c||_{p}$ to decay by a factor of $e^{-1}$. Based on the fit in Eq.~\eqref{eq:cnorm_fit}, we can calculate $T_{Pe}$ exactly as 
\begin{equation}
	T_{Pe} = \tau\Gamma(1+1/\alpha),
\label{eq:TPe_defn}
\end{equation}
where $\Gamma(z) := \int_0^\infty \zeta^{z-1}e^{-\zeta}\,\mathrm{d}\zeta$ is the Gamma function. The $Pe$ [$= L^2/(DT)$, recall Sect.~\ref{sec:diffusion}] subscript reminds us that this number depends on the relative ``strength'' of diffusion in the problem.
In the literature on finite Markov chains, the number of iterations $T_{Pe}$ would be called the \emph{stopping time}. 

\begin{figure}[t]
  \centering
  \includegraphics[width=0.5\textwidth]{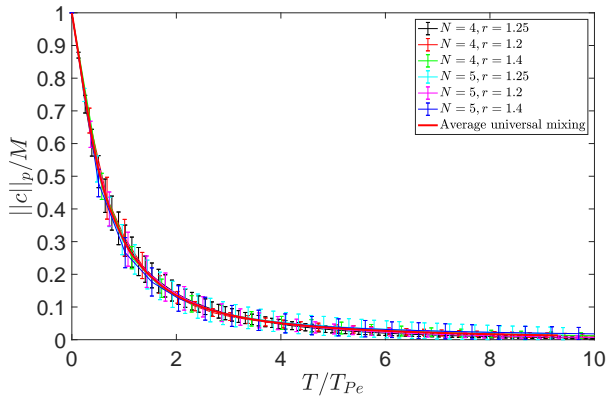}
  \caption{(Color online.) Rescaled mixing norm of the concentration/color versus the rescaled iterations, for a diffusion coefficient of $D=0.5$. All curves collapse as the horizontal axis is scaled by the stretched exponentials' e-folding time $T_{Pe}=\tau\Gamma(1+1/\alpha)$; $p=2$. The average universal mixing profile (bold curve) corresponds to the fit of the average of all the average rescaled mixing norm decay curves in this rescaled plot. The best-fit parameters, according to Eq.~\eqref{eq:cnorm_fit}, for the average curve are $\tau_\mathrm{universal} = 0.8706$ and $\alpha_\mathrm{universal} = 0.7920$.}
  \label{fig:half-mixing}
\end{figure}

To collapse all the mixing norm decay across different parameters, we now rescale each \emph{averaged} mixing norm decay curve (recall Sect.~\ref{sec:measures}) as $||c||_{p}(T) \mapsto ||c||_{p}(T)/M$ so that the mixing norm will range from $0$ to $1$; and, we rescale $T\mapsto T/T_{Pe}$. Then, we plot $||c||_{p}(T)/{M}$ versus $T/T_{Pe}$, for \emph{all} simulations that we have performed, onto the single plot shown in Fig.~\ref{fig:half-mixing}. A clear collapse of all the data is observed, when the horizontal axis is rescaled by $T_{Pe}=\tau\Gamma(1+1/\alpha)$, which suggest that  stopping time can be approximated by the e-folding time of the stretched exponential fit. In this plot, error bars denote one standard deviation from the mean (errors bars are bounded by $0$ from below, obviously) of all the concentration curves over all permutations that we have considered.  Error bars are used in order to be able to provide a sense for the behavior of all the possible ({\it i.e.}, across the allowed permutations $\Pi$) mixing norm decay profiles in a single plot. (Note that these are not error bars quantifying the uncertainty in the nonlinear least-squares fits.) 

In Fig.~\ref{fig:half-mixing}, each curve corresponds to a given IET with a fixed diffusion coefficient $D$, segment ratio $r$ and number of pieces $N$. The most enticing aspect of Fig.~\ref{fig:half-mixing} is that a \emph{single universal} profile of the form given in Eq.~\eqref{eq:cnorm_fit} can be fit to the average of all the average curves. It can be calculated that this ``average of averages'' profile, of $||c||_{p}(T)/{M}$ versus $T/T_{Pe}$, has  $\tau_\mathrm{universal} = 0.8706$, $\alpha_\mathrm{universal} = 0.7920$ (for the chosen $D=0.5$), which are now independent of $N$ and $r$ (unlike Table~\ref{tb:rLT} and Fig.~\ref{fig:different-average-norm})!

Though Fig.~\ref{fig:half-mixing} shows the relaxation of the line segment's color to the mean through \emph{many different} mixing protocols, it is clear that (after the appropriate rescaling) all IETs with diffusion behave in a \emph{universal} way. However, to provide further evidence of a cut-off phenomenon, we need to observe the transition from unmixed to mixed sharpening as $D$ becomes small (equivalently, as $Pe$ becomes large), as in \cite{lw08}. To do so, however, we need to estimate $T_{Pe}$ {\it a priori} on physical grounds, not by defining it via Eq.~\eqref{eq:TPe_defn}. When using the operational definition in Eq.~\eqref{eq:TPe_defn} to compute $T_{Pe}$, we find that the sensitivity to $Pe$ is weak, which is in line with the weak dependence observed for some maps in \cite{lw08}. Therefore, we would like to determine whether a physically-motivated prediction of $T_{Pe}$ allows us to more clearly see a sharpening of the concentration norm decay curves as $Pe\to\infty$ and, thus, to better substantiate the possibility of a cut-off.

\subsection{Predicting the stopping time $T_{Pe}$}
\label{sec:predict_TPe}

Schlick {\it et al.}~\cite{scuol13} proposed a simple one-dimensional analysis of diffusion between two subsegments of unequal color (on a normalized domain with a given $Pe$), using an analytic solution to the diffusion equation [{\it i.e.}, Eq.~\eqref{eq:diffusion_eq}]. (Muzzio and Ottino \cite{mo89b} previously considered the related case of reaction-diffusion.) In \cite{scuol13}, the line segment was taken to have length $2\ell$ ({\it i.e.}, each subsegment was of length $\ell$), and the colors were taken to be $c_1=0$ and $c_2=1$ without loss of generality. A solution was developed using eigenfunctions, from which it was determined that a subsegment length of $\ell^*$, where 
\begin{equation}
	\ell^* = \ell^*(\hat{T}) = \pi  \sqrt{\frac{\hat{T}}{2Pe}},
\label{eq:lstar}
\end{equation}
will be ``washed out'' in a characteristic normalized/dimensionless time $\hat{T}$ to be made precise below (given a specific P\'eclet number $Pe$) \footnote{Since Eq.~\eqref{eq:lstar} introduces a length $\ell^*\propto \sqrt{D}$, it can be considered as a type of \emph{Batchelor scale} \cite{b59,t08}, which describes the smallest length scale of fluctuations that can persist in a fluid flow before they are homogenized by molecular diffusion.}. Thus, based on this analysis from \cite{scuol13}, we pose the following question: when will the \emph{average} continuous-color subsegment length, denoted $\ell_m$, in our cutting and shuffling process \emph{without} diffusion reach $\ell^*$? This question is important because, if $\ell_m\simeq \ell^*$ then $\hat{T}$ iterations of the cutting and shuffling \emph{with} diffusion ($D\ne0$, for a given $Pe$) later, the concentration of the striation will be damped out (decrease) by $e^{-2}\approx 13.5\%$ \cite[p.~15]{scuol13}. In other words, $\hat{T}$ is the double-e-folding time of the advection--diffusion process.

For cutting and shuffling without diffusion ($D=0$), we first rescale the problem as in Sect.~\ref{sec:diffusion}, which yields the P\'eclet number definition (using a ``velocity'' scale $L/T_\mathrm{max}$):
\begin{equation}
	Pe = \frac{L^2}{DT_\mathrm{max}}. 
\end{equation}
Under this rescaling both $\hat{X} := X/L$ and $\hat{T} := T/T_\mathrm{max}$ run from $0$ to $1$. Then, the average subsegment length ({\it i.e.}, the average of the lengths of subsegments of continuous color) can be trivially shown to be given exactly by
\begin{equation}
	\ell_m (\hat{T}) = \frac{1}{\hat{C}(\hat{T})+1}.
\label{eq:lm}
\end{equation}
Here, as before, $\hat{C}(\hat{T}):=C(T)$ is the number cutting interfaces as defined after $\hat{T}$ normalized iterations so the number of distinct subsegment of continuous color is clearly $\hat{C}(\hat{T})+1$. Next, we seek to estimate the number of iterations $\hat{T}$ required for diffusion to ``wash out'' the color gradients.

\begin{figure}[t]
  \centering
  \includegraphics[width=0.5\textwidth]{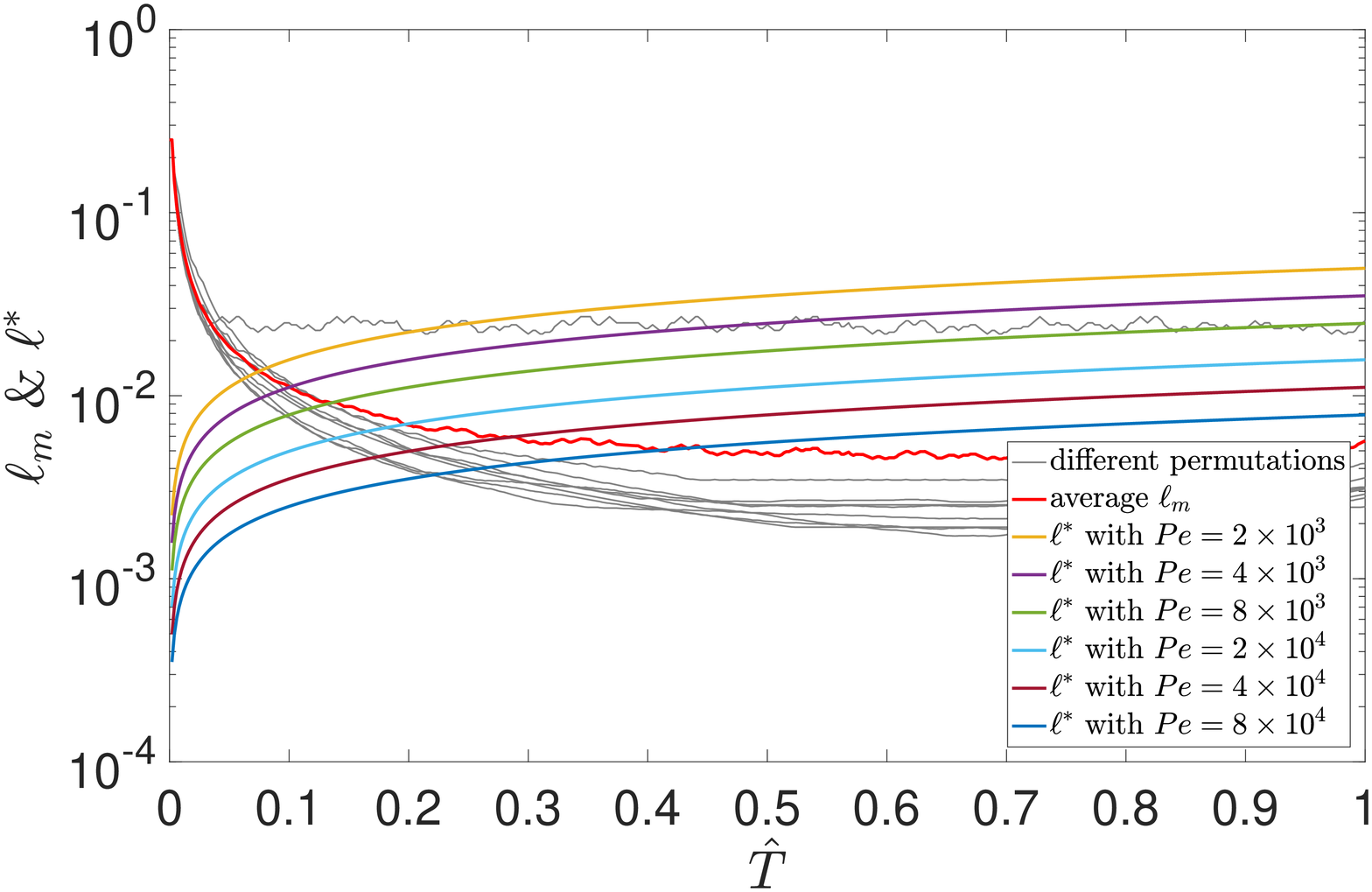}
  \caption{(Color online.) Graphical illustration of solving for $\hat{T}=\tilde{T}_{Pe}/T_\mathrm{max}$ such that $\ell^*(\hat{T}) = \ell_m(\hat{T})$. The average subsegment length $\ell_m$ (in the absence of diffusion) is shown as the light gray curves for $N=4$, $r=1.2$, $T_\mathrm{max} = 500$ and $D=0$ for various $\Pi$; the curve labeled ``average $\ell_m$'' is the average of the light gray curves ({\it i.e.}, over $\Pi$). The intersections of this averaged curve with the  $\ell^*$ curves (different colors correspond to different $Pe$ values, as in the legend) yield the values of $\tilde{T}_{Pe}/T_\mathrm{max} = 0.0740$, $0.102$, $0.136$, $0.196$, $0.284$, $0.408$.}
  \label{fig:l_star}
\end{figure}

To this end, in Fig.~\ref{fig:l_star}, we show visually how to determine when $\ell^*\simeq \ell_m$. In the absence of diffusion, the number of iterations required for the latter condition to hold is given by the $\hat{T}$ values at the intersections of the $\ell_m$ and $\ell^*$ curves in Fig.~\ref{fig:l_star}. At these $\hat{T}$, we can expect that diffusion dominates the dynamics, leveling the concentration gradients. Thus, we would like to argue that these values of $\hat{T}$ are  estimates of the stopping times. Let us introduce the notation $\hat{T}=\tilde{T}_{Pe}/T_\mathrm{max}$ (with the tilde introduced to clearly distinguish this value from the one in Sect.~\ref{sec:cutoff}) for this normalized stopping time, which is now defined based on Eqs.~\eqref{eq:lstar} and \eqref{eq:lm} as the solution of
\begin{equation}
\pi  \sqrt{\frac{\hat{T}}{2Pe}} = \frac{1}{\hat{C}(\hat{T})+1}.
\label{eq:TPe_new}
\end{equation}
Unfortunately, since the number of cutting interfaces $C(T)$ is a complicated function, for which we do not have a closed form solution, Eq.~\eqref{eq:TPe_new} must be solved numerically.

Notice that, it may turn out that for a given $Pe$, there is no solution to Eq.~\eqref{eq:TPe_new} (equivalently, the $\ell^*$ curve might never intersect the $\ell_m$ curve in Fig.~\ref{fig:l_star}). This situation occurs if the dynamics of the diffusionless IET are periodic, and the subsegment reassembles itself. In principle, our design rules given in Sect.~\ref{sec:choose_params} (specifically the exclusion of certain permutations), should preclude the possibility of periodic dynamics in $\ell_m$. If this were the case, nevertheless, then the average subsegment length never reaches the diffusion scale, thus for the given $Pe$ value, it is not expected that diffusion can significantly affect the mixing over the given $\hat{T}$ iterations [recall the discussion of Eq.~\eqref{eq:lstar}]. Indeed, we might expect that as $Pe\to\infty$, the solution to Eq.~\eqref{eq:TPe_new}, $\tilde{T}_{Pe}/T_\mathrm{max}\to\infty$ also.

\begin{figure}
  \centering
  \includegraphics[width=0.5\textwidth]{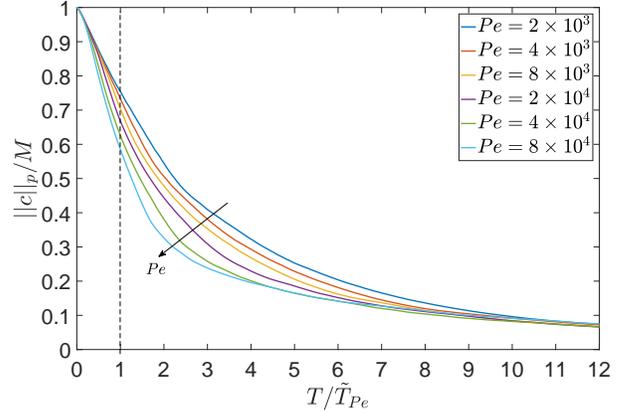}
  \caption{(Color online.) $Pe$ dependence of the average mixing norm $||c||_p(T)$ for $N=4$ and $r=1.2$, with $T$ rescaled using $\tilde{T}_{Pe}=45.01$, $58.00$, $79.00$, $124.0$, $142.0$, $189.09$ for $D = 0.451$, $0.225$, $0.113$, $0.0563$, $0.0375$, $0.0282$, respectively, as per Fig.~\ref{fig:l_star}.}
  \label{fig:new_rescale_TPe}
\end{figure}

After obtaining the value of $\tilde{T}_{Pe}/T_\mathrm{max}$ numerically from Eq.~\eqref{eq:TPe_new}, we can verify whether the foregoing argument about the influence of diffusion is valid by calculating the mixing norm decay with the specific $Pe$ given initially, from which we immediately get the corresponding diffusion coefficient to be used in a cutting and shuffling simulation:
\begin{equation}
	D = \frac{L^2}{Pe\, T_\mathrm{max}}.
\label{eq:D_new}
\end{equation}
To summarize: supposing a P\'eclet number (inverse dimensionless diffusivity) is known for a line segment of length $L$ normalized to $1$, then $\ell^*$ is estimated by Eq.~\eqref{eq:lstar} based on \cite{scuol13}. Next, the number of normalized iterations $\tilde{T}_{Pe}/T_\mathrm{max}$ of the diffusionless ($D=0$) cutting and shuffling process until diffusion would ``take over'' is estimated from Eq.~\eqref{eq:TPe_new}, from which $\tilde{T}_{Pe}$ trivially follows. Next, to verify the latter is an estimate of the stopping time, a cutting and shuffling simulation with diffusion is performed, using the properly matched diffusivity according to Eq.~\eqref{eq:D_new}.

\begin{figure}
  \centering
  \includegraphics[width=0.5\textwidth]{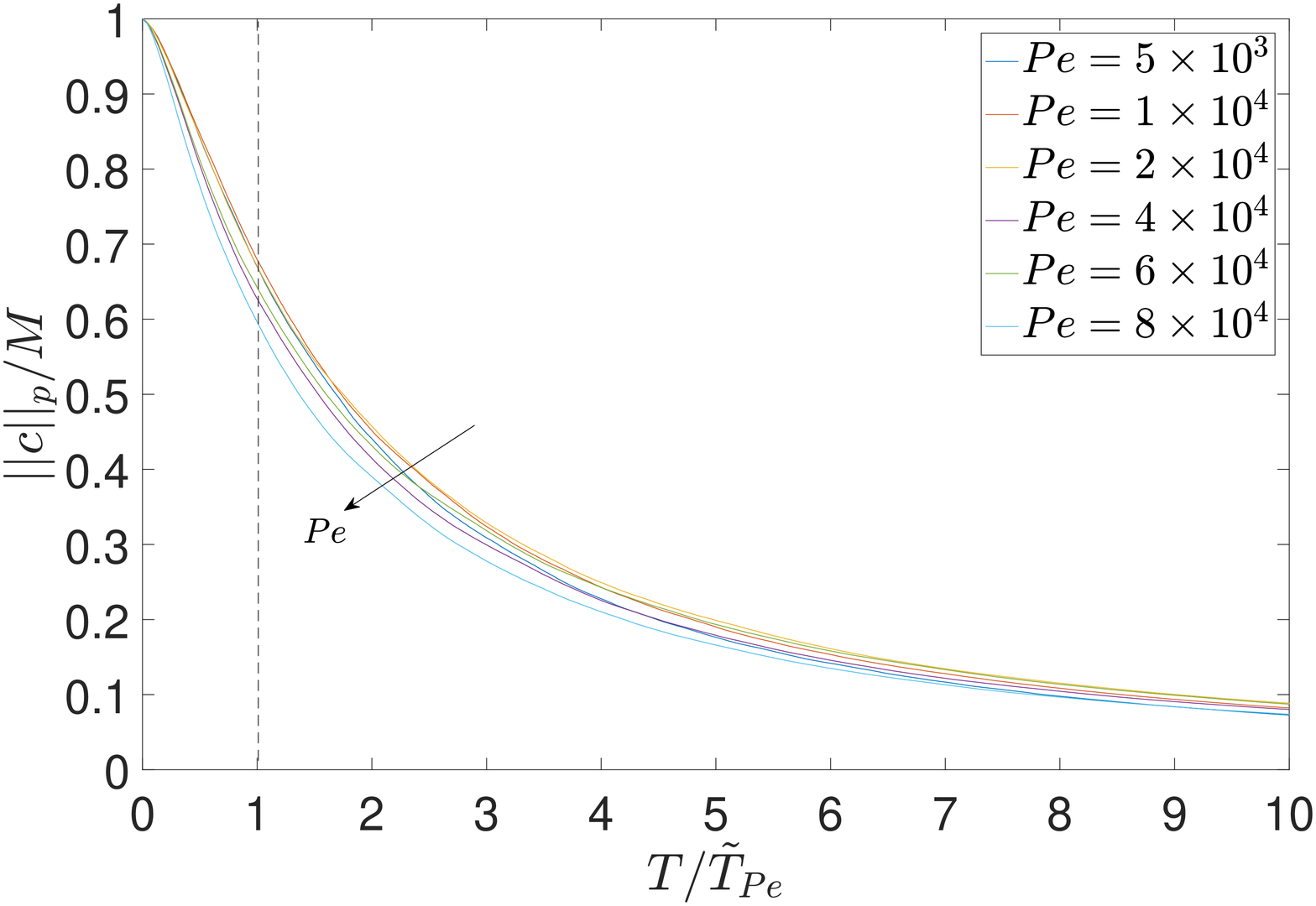}
  \caption{(Color online.) $Pe$ dependence of the average mixing norm $||c||_p(T)$ for $N=4$ and $r=1.4$, with $T$ rescaled using $\tilde{T}_{Pe}=37.0$, $51.0$, $68.0$, $98.0$, $142$, $204$ for $D = 0.45024$, $0.22512$, $0.11256$, $0.045024$, $0.022512$, $0.011256$, respectively.} 
  \label{fig:new_rescale_TPe2}
\end{figure}

\begin{figure}
  \centering
  \includegraphics[width=0.5\textwidth]{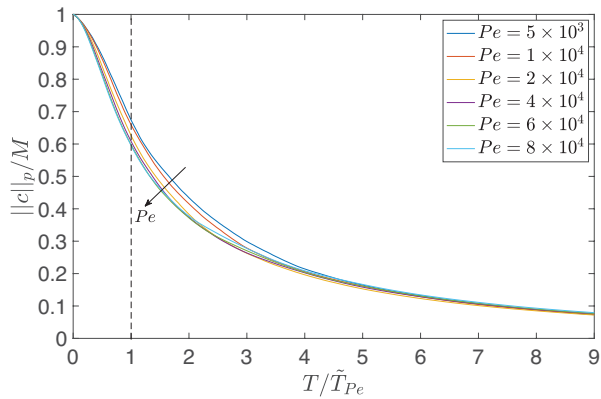}
  \caption{(Color online.) $Pe$ dependence of the average mixing norm $||c||_p(T)$ for $N=4$ and $r=1.6$, with $T$ rescaled using $\tilde{T}_{Pe}=49.0$, $65.0$, $91.0$, $124$, $151$, $172$ for $D = 0.446$, $0.223$, $0.112$, $0.0558$, $0.0372$, $0.0279$, respectively.}
  \label{fig:new_rescale_TPe3}
\end{figure}

A result from this numerical approach is illustrated in Fig.~\ref{fig:new_rescale_TPe} for a select choice of $N$ and $r$. Clearly, the otherwise disparate mixing norm curves are grouped together when plotting against $T/\tilde{T}_{Pe}$ for the estimated values of $\tilde{T}_{Pe}$ obtained from Eq.~\eqref{eq:TPe_new} and as shown visually in Fig.~\ref{fig:l_star}. More importantly, however, the grouping shows a ``steepening'' of the profiles as $Pe \to \infty$. This steepening (as in \cite{lw08}) is suggestive of a cut-off developing (a sharp transition from an unmixed state to a mixed state), the ideal form of which is represented by the dashed vertical line connecting $1$ to $0$ at $T/\tilde{T}_{Pe}=1$ in Fig.~\ref{fig:new_rescale_TPe}. 

Figures~\ref{fig:new_rescale_TPe2} and \ref{fig:new_rescale_TPe3} show a similar result, again for $N=4$ subsegments, but with ratios $r=1.4$ and $r = 1.6$, respectively; although less pronounced, the cut-off phenomenon appears to be present. The differences can be attributed to the improved mixing that occurs as $r$ becomes ``more irrational'' (increasing continued fraction expansion from $1.2$ to $1.4$ to $1.6$ as in Table~\ref{table:norm-fit-parameter}). Thus, the case of $r=1.2$ in Fig.~\ref{fig:new_rescale_TPe} is, in a sense, ``special'' in clearly showing the steepening. Nevertheless, a steepening with increasing $Pe$ is observed in all three  Figs.~\ref{fig:new_rescale_TPe}, \ref{fig:new_rescale_TPe2} and \ref{fig:new_rescale_TPe3}, providing numerical evidence suggestive of the cut-off phenomenon, across different choices of $N$ and $r$ in our cutting and shuffling process with diffusion.

\begin{figure}
  \centering
  \includegraphics[width=0.5\textwidth]{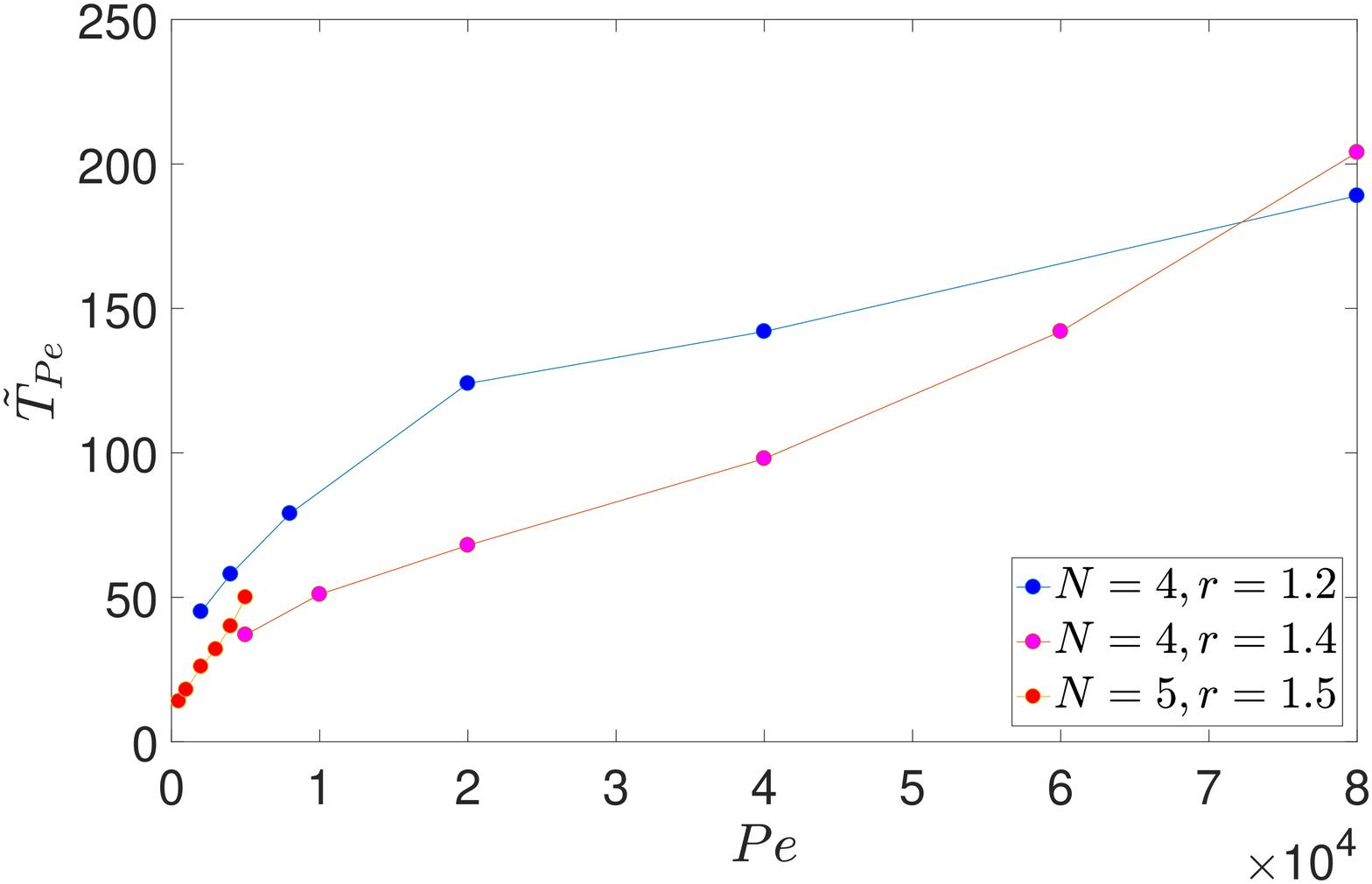}
  \caption{(Color online.) Dependence of the stopping time $\tilde{T}_{Pe}$, as calculated through Eq.~\eqref{eq:TPe_new}, as a function of $Pe$ for three combinations of $N$ and $r$.}
  \label{fig:tpe_pe}
\end{figure}

Finally, having considered the effect of different $N$ and $r$, we would like to compare the various predictions in the above discussion to understand the relationship between $\tilde{T}_{Pe}$ and $Pe$. This dependence is shown in  Fig.~\ref{fig:tpe_pe}. In general, the estimated stopping time $\tilde{T}_{Pe}$ increases with $Pe$, which should be expected given the hypothesis of a cut-off phenomenon. However, the specific dependence appears to be sensitive upon the choice of $N$ and $r$, and no clear pattern emerges. Thus, it remains the subject of future inquiry whether a specific functional form could be determined to specify the relationship between $\tilde{T}_{Pe}$ and $Pe$ {\it a priori}.

\section{Conclusion}
\label{sec:conclusion}

In the present work, we discussed the effects of incorporating diffusion into one-dimensional cutting and shuffling maps represented by interval exchange transformations. The most obvious conclusion is that diffusion leads to significantly enhanced mixing compared to cutting and shuffling alone. However, in the presence of diffusion, we must additionally be careful how we quantify mixing. Specifically, a ``mixing norm'' is a more effective way to quantify the degree of mixing compared to the number of cutting interfaces, a quantity of interested in the dynamical systems theory of interval exchange transformations.

Next, we proposed a parametrization of the possible mixing behaviors. Indeed, the class of cutting and shuffling protocols considered has a large parameter space: the number of initial subsegments $N$, the adjacent subsegment length ratio $r$, the shuffling permutation $\Pi$, and the diffusivity $D$ can all be varied independently. Our parametrization consists of fitting the decay of the mixing norm with the number of iterations to two parameters: a time constant $\tau$ and a skewness parameter $\alpha$, which were both found to depend on the ratio $r$ and on the diffusion coefficient $D$, for fixed $N$. Through this approach, we showed that, even though a large number of dynamical behaviors are possible, an appropriate rescaling of the mixing norm decay curves leads to a \emph{universal} mixing curve, describing (within some error margin) all cutting and shuffling protocols. This universality rests upon the fact that a number of iteration for the mixing norm to decay by a factor of $e^{-1}$, denoted $T_{Pe}$, can be found analytically as the e-folding time of the stretched exponential fit of the mixing norm decay curve for each protocol.

Another question we sought to address is whether cut-offs, and the concomitant concept of stopping times, from finite Markov chain theory apply here. To this end, we sought to determine a critical number of iterations $\tilde{T}_{Pe}$ (the stopping time) at which diffusion would ``kick in'' thus homogenizing the mixture. In doing so, we explored the limit of vanishing diffusivity ({\it i.e.}, $Pe\to\infty$), providing evidence (at least numerically) that the transition from an unmixed stated becomes sharper as $Pe\to\infty$. Our numerical exploration provides initial evidence that the concepts of cut-offs and stopping times, which elegantly explain that a deck of fifty-two cards requires about seven (and not any more) shuffles to become randomized, are relevant to interval exchange transformations with diffusion. Of course, unlike card shuffling, interval exchange transformations possess significant complexity. Even though such maps are easy to describe qualitatively, their mathematical theory remains an advanced and difficult topic in dynamical systems. Our results on cut-offs are also distinct from those by Liang and West \cite{lw08} because they considered \emph{chaotic} maps in several dimensions, while our interval exchanges are at best weakly mixing (though never truly so since we work on finite grids and with integer arithmetic).

A possible avenue of future work involves extending our cutting and shuffling approach with diffusion to consider chemical reactions occurring between the subsegments of different color. Such an extension could connect to the classical work on evolution and coarsening of lamellar structures in chaotic mixing \cite{mo89,mo89b,sb91,ccr99}. (Lamellar models of mixing remain an active topic of  research today \cite{szllbm18}.) The interplay between the lamellar width distribution, coupled to chemical reactions and diffusion processes, plays a key role in the evolution toward a steady state and, thus, the final yield of a chemical reaction. Clifford {\it et al.}~\cite{ccr99} discussed these issues at length, however, overall they initially considered only ``simple'' initial arrangement. Specifically, Clifford {\it et al.}\ noted that a weakness of their approach was that ``lamellae can have only two different widths, while realistic fluid flows generate lamellae with a wide range of widths'' \cite[p.~305]{ccr99}. Meanwhile, earlier work by Sokolov and Bluman concluded that ``the course of reaction is governed mainly by mixing and not by diffusion or kinetics.'' \cite[pp.~3698--3699]{sb91}. Thus, cutting and shuffling a line segment with reaction and diffusion presents a natural model in which to capture such complexity. 

As we have shown in the present work, when permutations leading to poor mixing and pathological cases are excluded, the behavior of the remaining protocols of cutting and shuffling with diffusion is, in a sense, universal although the ``stopping time'' $T_{Pe}$ is highly sensitive to the details of the protocol. Conceivably, such universality persists in the presence of reactions with the stopping time becoming the quantity one may wish to optimize. Indeed, the lamellar distributions under pure reaction and diffusion have been shown to evolve in a \emph{self-similar} manner \cite{mo89,mo89b}, suggesting some level of universality already exists in the process. Nevertheless, these remain questions that must be addressed in future work. 

Finally, we note in passing that a recently introduced thermodynamic approach by Brassart {\it et al.}~\cite{bls16,bls18} to molecular mixing makes use of the same basic concepts as in the continuum theories of the mixing of fluids: shear, dilation, diffusion and ``swap.'' Brassart {\it et al.}~\cite{bls16,bls18} define ``swap'' as ``preserv[ing] the shape and the volume, but chang[ing] the ratio of the two species of molecules'' in a piece of material \cite[p.~50]{bls16}. Indeed, this definition coincides with the notion of cutting and shuffling discussed herein if we consider how the color changes with iterations on an arbitrary finite piece of the original line segment.

\acknowledgments
I.C.C.\ would like to thank Jean-Luc Thiffeault for helpful discussions on cut-offs in Markov chains and for pointing out the work by Liang and West \cite{lw08} to him. The authors also thank Lachlan Smith and the anonymous referees for helpful comments and criticisms, which improved the manuscript. Specifically, we thank one of the anonymous referees for pointing out the relevance of the e-folding time of the stretched exponential fit.

\bibliography{c_s_refs}

\end{document}